\documentclass{gtart_h}  

\def\ifplaintex{\expandafter\ifx\csname documentclass\endcsname\relax}

\ifplaintex 
\else
\fi

\expandafter\ifx\csname beginpicture\endcsname\relax
\expandafter\ifx\csname documentclass\endcsname\relax
\input pictex \else
\input prepictex \input pictex \input postpictex \fi\fi

\def\gt{{\mathsurround=0pt\it $\cal G\mskip-2mu$eometry \&\ 
$\cal T\!\!$opology}}        %  journal title in recommended style

\def\gtp{{\mathsurround=0pt\it $\cal G\mskip-2mu$eometry \&\ 
$\cal T\!\!$opology $\cal P\!$ublications}}  % GT publications

\def\lognumber#1{\def\thelognumber{#1}}
\def\volumenumber#1{\def\thevolumenumber{#1}}
\def\papernumber#1{\def\thepapernumber{#1}}
\def\volumeyear#1{\def\thevolumeyear{#1}}

\def\pagenumbers#1#2{\def\startpage{#1}\def\finishpage{#2}}
\def\published#1{\def\publishdate{#1}}
\def\proposed#1{\def\theproposer{#1}}
\def\seconded#1{\def\theseconders{#1}}
\def\received#1{\def\receiveddate{#1}}

\def\accepted#1{\def\accepteddate{#1}}

\def\asciiaddress#1{\def\theasciiaddress{#1}}
\def\asciiemail#1{\def\theasciiemail{#1}}

\def\asciikeywords#1{\def\theasciikeywords{#1}}

\let\\\par\let\thelognumber\relax
\let\thevolumenumber\relax\let\thepapernumber\relax
\let\thevolumeyear\relax\let\thesamplenumber\relax\let\startpage\relax
\let\finishpage\relax\let\publishdate\relax\let\receiveddate\relax
\let\reviseddate\relax\let\accepteddate\relax\let\theasciititle\relax
\let\theasciiauthors\relax\let\theasciiaddress\relax
\let\theasciiabstract\relax\let\theasciikeywords\relax
\let\theasciiemail\relax\let\theshortauthors\relax\let\theshorttitle\relax

\long\def\maketitlep{   % start of definition of \maketitlep

\count0=\startpage

\gt\hfill      %   Journal title (top left) 
\beginpicture
\setcoordinatesystem units <0.33truein, 0.33truein> point at 2.2 0.9
\setplotsymbol ({$\cal G$})
\plotsymbolspacing=9truept
\circulararc 315 degrees from 0 1 center at 0 0
\setplotsymbol ({$\cal T$})
\circulararc 315 degrees from 1 -1 center at 1 0
\endpicture
\break
{\small\ifx\thesamplenumber\relax % sample?  
Volume \else Sample
\fi\thevolumenumber\ (\thevolumeyear)
\startpage--\finishpage\nl
Published: \publishdate}
\vglue 0.5truein plus 0.4fil minus 0.1truein

{\parskip=0pt\leftskip 0pt plus 1fil\def\\{\par\smallskip}{\ifplaintex\large
\else\Large\fi\bf\thetitle}\par\medskip}   

\vglue 0pt plus 0.1fil 

{\parskip=0pt\leftskip 0pt plus 1fil\def\\{\par}{\sc\theauthors}
\par\medskip}

\vglue 0pt plus 0.1fil 

{\small\parskip=0pt\let\newline\\
{\leftskip 0pt plus 1fil\def\\{\par}{\sl\theaddress}\par}
\expandafter\ifx\theemail\relax    % email address?
\relax\else\vglue 5pt plus 0.02fil minus 2pt\def\\{\stdspace{\rm 
and}\stdspace} 
\cl{Email:\stdspace\tt\theemail}\fi
\ifx\theurl\relax                  % URL given?
\relax\else\vglue 5pt plus 0.02fil minus 2pt\def\\{\stdspace{\rm 
and}\stdspace}
\cl{URL:\stdspace\tt\theurl}\fi\par}

\vglue 7pt plus 0.3fil minus 3pt

{\bf Abstract}
\vglue 5pt plus 0.1fil minus 2pt

\theabstract

\vglue 7pt plus 0.3fil minus 3pt

{\bf AMS Classification numbers}\quad Primary:\quad \theprimaryclass

Secondary:\quad \thesecondaryclass

\vglue 5pt plus 0.3fil minus 2pt

{\bf Keywords:}\quad \thekeywords

\vglue 10pt plus 0.5fil minus 5pt

{\small  Proposed: \theproposer\hfill Received: \receiveddate\nl
Seconded: \theseconders\hfill 
\ifx\reviseddate\relax                         % paper revised?
Accepted: \accepteddate                        % no
\else
Revised: \reviseddate                          % yes
\fi}
\eject
}       %  end of definition of \maketitlep

\let\maketitlepage\maketitlep
\let\maketitle\maketitlepage

\font\phead=cmsl9 scaled 950
\font=cmsl9 scaled 1050
\font\pnum=cmbx10 scaled 913
\font\lnum=cmbx10 
\font\pfoot=cmsl9 scaled 950
\font=cmsl9 scaled 1050
\ifplaintex
\headline{\vbox to 0pt{\vskip -4.5mm\line{\small\phead\ifnum
\count0=\startpage ISSN 1364-0380 (on line)
1465-3060 (printed) \hfill {\pnum\folio}\else\ifodd\count0\def\\{ }% 
\ifx\theshorttitle\relax\thetitle\else\theshorttitle\fi\hfill{\pnum\folio}
\else\def\\{ and }{\pnum\folio}\hfill\ifx\theshortauthors\relax\theauthors
\else\theshortauthors\fi\fi\fi}\vss}}
\footline{\vbox to 0pt{\vglue 0mm\line{\small\pfoot\ifnum\count0=\startpage
\copyright\ \gtp\hfill\else
\gt, Volume \thevolumenumber\ (\thevolumeyear)\hfill\fi}\vss
}}
\else
\makeatletter
\def\@oddhead{{\small\lhead\ifnum\count0=\startpage ISSN 1364-0380 (on line)
1465-3060 (printed) \hfill {\lnum\number\count0}\else\ifodd\count0
\def\\{ }\ifx\theshorttitle\relax \thetitle \else\theshorttitle\fi\hfill
{\lnum\number\count0}\else\def\\{ and }{\lnum\number\count0}
\hfill\ifx\theshortauthors\relax 
\theauthors\else\theshortauthors\fi\fi\fi}}\def\@evenhead{\@oddhead}
\def\@oddfoot{\small\lfoot\ifnum\count0=\startpage\copyright\ \gtp\hfill\else
\gt, Volume \thevolumenumber\ (\thevolumeyear)\hfill\fi}
\def\@evenfoot{\@oddfoot}
\makeatother
\fi

\newwrite\gtoutfile
\long\gdef\makeheadfile{  %%% start of definition of \makeheadfile
{\def\\{, }\def\s{ }
\immediate\openout\gtoutfile head.xxx
\immediate\write\gtoutfile{Proxy-for: \ifx\theasciiauthors\relax
\theauthors\else\theasciiauthors\fi\s<\ifx\theasciiemail\relax\theemail\else\theasciiemail\fi>}
\immediate\write\gtoutfile{\noexpand\\}
\immediate\write\gtoutfile{Authors: \ifx\theasciiauthors\relax
\theauthors\else\theasciiauthors\fi}
{\def\\{ }\immediate\write\gtoutfile{Title: \ifx\theasciititle\relax
\thetitle\else\theasciititle\fi}}
\immediate\write\gtoutfile{Subj-class: GT or SG or MG etc}
\immediate\write\gtoutfile{MSC-class: \theprimaryclass\ifx\thesecondaryclass\relax\else, \thesecondaryclass\fi}
\immediate\write\gtoutfile{Journal-ref: Geom. Topol. \thevolumenumber
(\thevolumeyear) \startpage-\finishpage}
\immediate\write\gtoutfile{Comments: Published by Geometry and Topology at}
\immediate\write\gtoutfile{\s\s http://www.maths.warwick.ac.uk/gt/GTVol\thevolumenumber/paper\thepapernumber.abs.html}
\immediate\write\gtoutfile{\noexpand\\}
\immediate\write\gtoutfile{}
\ifx\theasciiabstract\relax
\immediate\write\gtoutfile{\theabstract}\else
\immediate\write\gtoutfile{\theasciiabstract}\fi
\immediate\write\gtoutfile{}
\immediate\write\gtoutfile{\noexpand\\}
\immediate\write\gtoutfile{}
\immediate\closeout\gtoutfile}}  %%% end of definition of \makeheadfile

\def\maketitlepage{\maketitlep\makeheadfile}
\let\maketitle\maketitlepage

\lognumber{490}
\received{31 August 2004}
\volumenumber{9}\papernumber{12}\volumeyear{2005}
\pagenumbers{375}{402}   
\published{3 March 2005}
\accepted{28 February 2005}
\proposed{Walter Neumann}
\seconded{Wolfgang Metzler, Cameron Gordon}

\usepackage{amsmath,amssymb}

\newtheorem{thm}{Theorem}[section]
\newtheorem{lem}[thm]{Lemma}
\newtheorem{prop}[thm]{Proposition}

\newtheorem{theor}{Theorem}

\theoremstyle{remark}
\newtheorem{notation}[thm]{Notation}
\newtheorem{defn}[thm]{Definition}
\newtheorem{conv}[thm]{Convention}
\newtheorem{rem}[thm]{Remark}

\renewcommand{\rk}{\operatorname{rk}}
\newcommand{\ab}{\operatorname{ab}}
\newcommand{\Comm}{\operatorname{Comm}}

\begin{document}

\title{Kleinian groups and the rank problem}

\authors{Ilya Kapovich\\Richard Weidmann}

\address{Department of Mathematics, University of Illinois at
  Urbana-Champaign\\1409 West Green Street, Urbana, IL 61801, USA}
\secondaddress{Fachbereich Mathematik, Johann Wolfgang Goethe
  Universit\"at\\Robert Mayer-Stra\ss e 6--8, 60325 Frankfurt, Germany}

\asciiaddress{Department of Mathematics, University of Illinois at
  Urbana-Champaign\\1409 West Green Street, Urbana, 
IL 61801, USA\\and\\Fachbereich Mathematik, Johann Wolfgang Goethe
  Universitat\\Robert Mayer-Strasse 6-8, 60325 Frankfurt, Germany}

\gtemail{\mailto{kapovich@math.uiuc.edu}{\rm\qua 
and\qua}\mailto{rweidman@math.uni-frankfurt.de}}
\asciiemail{kapovich@math.uiuc.edu, rweidman@math.uni-frankfurt.de}


\asciikeywords{Word-hyperbolic groups, Nielsen methods, 3-manifolds}

\begin{abstract}
  We prove that the rank problem is decidable in the
  class of torsion-free word-hyperbolic Kleinian groups. We also show
  that every group in this class has only finitely many Nielsen
  equivalence classes of generating sets of a given cardinality.
\end{abstract}

\maketitlepage

\section{Introduction}

If $G$ is a finitely generated group, the \emph{rank of $G$}, denoted
$\rk(G)$, is the smallest cardinality of a subset $S\subseteq G$ such
that $S$ generates $G$. The \emph{rank problem} for a particular class
of finitely presented groups asks if there is an algorithm that, given
a finite presentation of a group $G$ from the class, computes the rank
of~$G$.

The rank problem is one of the more enigmatic and little understood
group-theoretic decision problems. For example, the word problem, the
conjugacy problem and even (in the torsion-free case) the isomorphism
problem~\cite{Se} are solvable for the class of word-hyperbolic
groups. On the other hand, by a theorem of Baumslag, Miller and
Short~\cite{BMS}, the rank problem is unsolvable for word-hyperbolic
groups. The main ingredient in their proof is a remarkable
construction due to Rips~\cite{R}.  Namely that, given an
arbitrary finitely presented group $Q$, there exists a short exact
sequence
\[1\to K\to G\to Q\to 1\] such that $G$ is torsion-free, non-elementary,
word-hyperbolic and such that $K$ is two-generated. Using the
classical undecidability results, it is possible to create a family of
finitely presented groups where it is undecidable whether $Q$ is
trivial or not, and where every nontrivial $Q$ has rank at least
three. Then the group $G$ based on $Q$ via the Rips construction can
be generated by two elements if and only if $Q$ is trivial. Hence it
is undecidable whether $\rk(G)\le 2$ and therefore $\rk(G)$ is not
computable.

A crucial feature of Rips' construction is that if both $K$ and $Q$
are infinite then $K$ is not quasiconvex in $G$. As it turns out, it
is the presence of finitely-generated non-quasiconvex subgroups that
is often responsible for undecidability of various algorithmic
problems in the context of word-hyperbolic groups. Thus Kapovich and
Weidmann~\cite{KW04} proved that the rank problem is solvable for
torsion-free locally quasiconvex word-hyperbolic groups.

However, some of the most important and interesting examples of
word-hyper\-bolic groups come from the world of Kleinian groups and they
are not necessarily locally quasiconvex. For instance, let $M$ be a
closed hyperbolic 3--manifold fibering over a circle. In this case the
fiber is a closed surface of negative Euler characteristic.  Then
$G=\pi_1(M)$ splits as a semi-direct product $G=H\ltimes \mathbb Z$,
where $H$ is the fundamental group of a fiber. Thus $H$ is not
quasiconvex in $G$.

The rank problem for 3--manifold groups is particularly interesting
because of the connection between the rank and another important
invariant, the Heegaard genus. Thus if $M$ is a closed 3--manifold,
$G=\pi_1(M)$ and $h(M)$ is the Heegaard genus of $M$ then $\rk(G)\le
h(M)$. Waldhausen asked if in fact $\rk(G)$ is always equal to $h(M)$.
Boileau and Zieschang~\cite{BZ} constructed a counter-example by
producing a family of Seifert manifolds of Heegaard genus three with
2--generated fundamental groups. To this day the Waldhausen question
remains open in the case of hyperbolic manifolds (see, for example, a
paper of Dunfield and Thurston~\cite{DT} where some experimental data
is discussed).

The genus problem for 3--manifolds has been solved by
Johannson~\cite{J} for sufficiently large 3--manifolds. Moreover he
provides a procedure that produces all Heegaard splittings of a given
genus. The case of small 3--manifolds has been claimed by Jaco and
Rubinstein \cite{JR,Ru} and relies on the existence and
constructibility of so-called 1--efficient triangulations of
irreducible atoroidal 3--manifolds.

\begin{defn}
  Let $\mathcal M$ be the class of all torsion-free word-hyperbolic
  groups $G$ such that $G$ admits a properly discontinuous isometric
  action on $\mathbb H^3$.
\end{defn}

Our main result is:

\begin{theor}\label{main}
  There exists an algorithm which, given a finite presentation of a
  group $G$ from the class $\mathcal M$, finds the rank of $G$.
\end{theor}

The main technical tool needed for the proof of Theorem~\ref{main} is
machinery developed by Kapovich and Weidmann in \cite{KW,KW04} that
provides a far-reaching generalization of Nielsen's methods in the
general context of groups acting by isometries on Gromov-hyperbolic
spaces. Another important ingredient is the ``tameness conjecture'' for
Kleinian groups that has been recently proved by Agol~\cite{A} and,
independently, by Calegari and Gabai~\cite{CG}. Together with the
results of Canary~\cite{C} this yields a precise characterization of
non-quasiconvex subgroups in the Kleinian groups context.

\smallskip There are very few other classes of groups where the rank
problem is known to be solvable. These include, in particular, the
class of finitely generated nilpotent groups. For an arbitrary group
$G$ it is easy to see that $\rk(G)=\rk(G/\mathcal F(G))$, where
$\mathcal F(G)$ is the Frattini subgroup of $G$.  If $G$ is finitely
generated nilpotent, then the Frattini subgroup $\mathcal F(G)$
contains the commutator $[G,G]$.  This implies that $\rk(G)=\rk(G_{\ab})$
where $G_{\ab}=G/[G,G]$ is the abelianization of $G$. The rank problem
is also decidable for torsion-free word-hyperbolic locally quasiconvex
groups (Kapovich--Weidmann~\cite{KW04}), for (sufficiently large)
Fuchsian groups~\cite{W04} and for one-relator groups~\cite{KMS}. The
results of Arzhantseva and Ol'shanskii~\cite{AO} show that, in a
sense, the rank problem is ``generically'' solvable for finitely
presented groups. The work of Weidmann~\cite{W04} (see also
Kapovich--Weidmann~\cite{KW04}) shows that the presence of torsion
often creates a substantial difficulty for solving the rank problem.
Indeed, there are many natural and seemingly easy to understand
classes of finitely generated groups where the rank problem remains
open. These include: virtually abelian groups, virtually free groups,
virtually nilpotent groups, 3--manifold groups as well as lattices in
$\mathbb H^3$.

We should stress that Theorem~\ref{main} is an abstract computability
result. The nature of the proof is such that it cannot provide any
complexity bound on the running time of the algorithm. In particular,
the proof involves many ``general enumeration'' arguments, where several
procedures are run in parallel, and at least one of them is guaranteed
to eventually terminate, but where no complexity estimate is possible.
We also stress that all the algorithms considered in this paper need
to be \emph{uniform} in all of their arguments, including the group
presentation.

We also obtain:

\begin{theor}\label{B}
  Let $G$ be a group from the class $\mathcal M$. Then for every $k\ge
  1$ there are only finitely many Nielsen-equivalence classes of
  $k$--tuples generating $G$.
\end{theor}

A similar but stronger statement was obtained by Kapovich and
Weidmann~\cite{KW04} for torsion-free one-ended locally quasiconvex
hyperbolic groups. Namely, for such a group $G$ for every $k\ge 1$
there are, up to conjugacy, only finitely many Nielsen-equivalence
classes of $k$--tuples generating one-ended subgroups. This stronger
statement is false in the presence of non-quasiconvex subgroups and in
particular it fails for some groups from the class $\mathcal M$.
Indeed, let $M$ be a closed hyperbolic 3--manifold fibering over a
circle. Then $G=\pi_1(M)$ splits as a semi-direct product $G=H\ltimes
\langle t\rangle$ where $H$ is the fundamental group of a fiber. If
$k$ is the rank of $H$, then the subgroups $H_n:=\langle H,t^n\rangle$
are all $(k{+}1)$--generated. However, each $H_n$ is normal of index $n$
in $G$ and hence $H_n$ is not conjugate to $H_m$ in $G$ for $m\ne n$.
Thus there are, up to conjugacy, infinitely many Nielsen-equivalence
classes of $(k{+}1)$--tuples generating subgroups of finite index in $G$.

More recently Souto~\cite{So} has obtained an interesting result that
is relevant to the present article.  Namely, he proved in \cite{So}
that if $\phi$ is a pseudo-Anosov homeomorphism of a closed oriented
surface of genus $g\ge 2$, then the ranks of the fundamental groups of
the mapping tori of sufficiently high powers of $\phi$ are equal to
$2g+1$ and that for these groups there is only one Nielsen equivalence
class of generating $(2g{+}1)$--tuples. His proof uses ideas closely
related to those used in this article.

The proof of Proposition~\ref{vc} about virtually cyclic groups is
based on an argument explained to the first author by Derek Holt and
we express special thanks to him. We are also grateful to Brian
Bowditch, Misha Kapovich and Walter Neumann for helpful discussions
about 3--manifold groups.

The first author was supported by the Max Planck Institute
of Mathematics in Bonn and by NSF grant DMS--0404991.

\section{Algorithms in hyperbolic groups}

If $G$ is a group with a finite generating set $S$, we will denote by
$X(G,S)$ the Cayley graph of $G$ with respect to $S$.

We refer the reader to \cite{ABC,E-T,GH,GS,G,KS} for the basic
background information regarding word-hyperbolic groups,
quasiconvexity, Gromov hyperbolic spaces and their boundaries. Recall,
however, that a subgroup $H$ of a word-hyperbolic group $G$ is
\emph{quasiconvex} in $G$ if for some (equivalently, for any) finite
generating set $S$ of $G$ the subgroup $H$ is a quasiconvex subset of
$X(G,S)$, that is, there is $\epsilon\ge 0$ such that for any
$h_1,h_2\in H$ any geodesic $[h_1,h_2]$ in $X(G,S)$ is contained in
the $\epsilon$--neighborhood of $H$.

In subsequent sections we will be relying on the
following statement summarizing some known general algorithmic results
regarding hyperbolic groups. Proposition~\ref{summ} will be often used
implicitly throughout the paper.

\begin{prop}\label{summ}
  The following statements hold:

\begin{enumerate}
  
\item There is an algorithm that, given a finite presentation
  $G=\langle S|R\rangle$ of a word-hyperbolic group $G$, computes an
  integer $\delta\ge 0$ such that the Cayley graph $X=X(G,S)$ of $G$
  with respect to $S$ is $\delta$--hyperbolic. Moreover:

  {\rm(a)}\qua The algorithm then computes a finite state automaton accepting
  the language $L=L(G,S)$ of all short-lex geodesic words over $S$ for
  $G$. The algorithm then decides if $G$ is finite or infinite, and if
  $G$ is finite, the algorithm computes the order of $G$ and decides
  whether or not $G$ is virtually cyclic.
  
  {\rm(b)}\qua Given an arbitrary word $w$ over $S$, the algorithm decides
  whether or not $w=1$ in $G$, and computes the order of the element
  of $G$ represented by~$w$.

\item There is an algorithm with the following property.
  
  {\rm(a)}\qua Given a finite presentation $G=\langle S|R\rangle$ of a
  word-hyperbolic group $G$ and a finite set $Q$ of words over $S$
  generating a subgroup $H=\langle Q\rangle\le G$, the algorithm
  computes an integer $C\ge 0$ such that $H$ is a $C$-quasiconvex
  subset of $X(G,S)$ if $H$ is quasiconvex in $G$, and runs forever if
  $H$ is not quasiconvex in $G$.
  
  {\rm(b)}\qua If $H$ is quasiconvex in $G$, the algorithm also computes a
  finite state automaton $L_H=L_H(G,S)$ accepting all the short-lex
  geodesic words over $S$ representing elements of $H$.
  
  {\rm(c)}\qua If $H$ is quasiconvex in $G$, the algorithm decides whether
  $H=G$, that is to say, whether $Q$ generates $G$.
  
  {\rm(d)}\qua If $H$ is quasiconvex in $G$, the algorithm computes the index
  of $H$ in~$G$.
  
  {\rm(e)}\qua If $H$ is quasiconvex in $G$, the algorithm computes a finite
  presentation of $H$ on $S$.
  
  {\rm(f)}\qua If $Q_1$ is another set of words over $S$ generating a
  quasiconvex subgroup $H_1\le G$ (so that $H_1\cap H$ is also
  quasiconvex), the algorithm computes a finite generating set for
  $H_1\cap H$.
  
  {\rm(g)}\qua If $H\le G$ turns out to be quasiconvex, then, given an
  arbitrary word $w$ over $S$, the algorithm decides whether or not
  $w$ represents an element of~$H$.

\end{enumerate}

\end{prop}

\begin{proof}
  The algorithm detecting hyperbolicity of a finite presentation
  $G=\langle S|R\rangle$ and producing a hyperbolicity constant
  $\delta$ is due to Papasoglu~\cite{P}. A general result of
  \cite{E-T} states that, given a finite group presentation $G=\langle
  S|R\rangle$ known to possess a short-lex automatic structure with
  respect to $S$, one can algorithmically find such a structure. Since
  a word-hyperbolic group is short-lex automatic with respect to every
  finite generating set, this algorithm will always terminate if $G$
  is word-hyperbolic.

  The algorithm detecting quasiconvexity of a finitely generated
  subgroup is due to Kapovich~\cite{Ka96}. To decide if a hyperbolic
  group $G$ is virtually cyclic, we compute an automaton accepting the
  language $L(G,S)$ off all short-lex geodesic words over $S$ and
  check whether or not $L(G,S)$ has linear growth.

  All of the other statements, except (2d), follow from basic
  well-known facts about word-hyperbolic groups.
  
  To show that (2d) holds, suppose that $H$ is a quasiconvex subgroup
  of $G$. We may assume that $G$ is infinite since otherwise the
  problem is easily decidable.  We need to compute the index of $H$ in
  $G$.  Recall that all cyclic subgroups in a word-hyperbolic group
  are quasiconvex.
  
  For a quasiconvex subgroup $H$ of a word-hyperbolic group $G$ a
  result of Ar\-zh\-ant\-seva~\cite{Ar} implies that $[G:H]=\infty$ if and
  only if there is an element $g\in G$ of infinite order such that
  $\langle g\rangle \cap H=1$. Therefore we will run in parallel the
  following two procedures.
  
  We will start the Todd--Coxeter coset enumeration process~\cite{Si}
  for $H$ in $G$. If $[G:H]<\infty$, the process will eventually
  terminate and its output can be used to compute the index of $H$ in
  $G$.
  
  In parallel, we will start enumerating all elements $g$ of $G$ (that
  is all words over $S$). For each of them we check if $g$ has finite
  order or not. If $g$ has infinite order, compute the automaton
  accepting $L_{\langle g\rangle}$ and then an automaton accepting
  $L_{\langle g\rangle\cap H}=L_{\langle g\rangle}\cap L_H$. We then
  check whether or not $L_{\langle g\rangle\cap H}$ consists of only
  the empty word. If yes, we conclude that $\langle g\rangle \cap H=1$
  and hence $[G:H]=\infty$. If not, we go on to the next~$g$.
  
  Eventually one of these procedures will terminate and we will have
  computed the index of $H$ in $G$.
\end{proof}

We will further need the following simple lemma which applies in
particular to the class of word-hyperbolic groups.

\begin{lem}\label{enum:fi}
  Suppose $\mathcal P$ is a class of finitely presented groups with
  uniformly solvable word problem. There exists a partial algorithm
  with the following property.

  Let $G$ be a group from $\mathcal P$ given by a finite presentation
  $G=\langle S|R\rangle$.  Let $A\subset G$ be a finite subset of $G$
  generating a subgroup $H$.  Given a finite subset $Q\subseteq G$,
  the algorithm will eventually terminate if $H$ has finite index in
  $L:=\langle A,Q\rangle$ and it will run forever otherwise.

  Suppose further that $H$ turns out to have finite index in $L$, that
  $H$ is given by a finite presentation on $A$ and that the pair
  $(H,G)$ comes from a class of subgroups of groups from $\mathcal P$
  where the membership problem for $H$ in $G$ is solvable.  Then the
  algorithm will also produce a finite presentation of $L$ on $A\cup
  Q$.
\end{lem}

\begin{proof}
  Denote $B=A\cup Q$. Start enumerating all elements of $H$ (that is
  all words in $A$) and, in parallel, start enumerating all finite
  prefix-closed sets $W$ of words in $B$. For each such set
  $W=\{w_1,w_2,\dots, w_n\}$, using this enumeration of $H$ and the
  algorithm solving the word problem in $G$, start checking if it is
  true that for every $w_i$ and every $b\in B^{\pm 1}$ there is some
  $w_j$ such that $w_ibw_j^{-1}$ belongs to $H$, i.e. if $w_ib\in
  w_jH$. If yes, then, clearly, $H$ has finite index in $L$ as the set
  $W$ contains a representative of every coset of $H$ in $L$, i.e. as
  $L=\cup_{w\in W} Hw$. Conversely, it is also obvious that if $H$ has
  finite index in $L$, the algorithm will eventually discover it and
  terminate.
  
  Suppose now that $H$ was given by a finite presentation, that $G$
  has solvable membership problem with respect to $H$ and that $L$
  turns out to contain $H$ as a subgroup of finite index. Recall that
  the above algorithm produces a finite set $W$ such that
  $L=\cup_{w\in W} Hw$. Using the algorithm solving the membership
  problem for $H$ in $G$, we can decide when $Hw=Hw'$ and when $Hw\ne
  Hw'$ for all $w,w'\in W$. Using this information we can construct
  the Schreier coset graph for $H$ in $L$ with respect to the
  generating set $B$ of $L$. Combined with a finite presentation of
  $H$ this yields a finite presentation for $L$.
\end{proof}

\section{Reduction to the one-ended case}

In this section we will show that it suffices to proof the main
theorem for one-ended groups.

\begin{defn}
  Let $\mathcal M_1$ be the class of one-ended torsion-free groups $G$
  such that $G$ admits a properly discontinuous convex-cocompact
  isometric action on $\mathbb H^3$.
\end{defn}
It is clear that every group from $\mathcal M_1$ is word-hyperbolic so
that $\mathcal M_1\subseteq \mathcal M$. Also, each of the classes
$\mathcal M$ and $\mathcal M_1$ is closed under taking subgroups of
finite index.

Moreover, it is a straightforward corollary of Thurston's
Hyperbolization Theorem (see, for example, 
\cite[Theorem 1.43]{Ka01}) that every one-ended group from $\mathcal M$ actually
belongs to $\mathcal M_1$:

\begin{prop}\label{thur}
  Let $G$ be a group acting isometrically and properly discontinuously
  on $\mathbb H^3$ and such that $G$ is torsion-free word-hyperbolic
  and one-ended. Then $G$ is isomorphic to the fundamental group of a
  compact hyperbolic 3--manifold with (possibly empty) convex boundary.
  That is, $G$ belongs to $\mathcal M_1$.
\end{prop}

The above fact reduces the rank problem for $\mathcal M$ to the rank
problem for $\mathcal M_1$ because of Grushko's theorem, since by a
result of Gerasimov~\cite{Ger} in the class of word-hyperbolic groups
one can algorithmically compute a maximal free product decomposition
of a group into freely indecomposable factors.

\begin{prop}\label{ger}
  There exists an algorithm with the following properties. Given a
  finite presentation $G=\langle S|R\rangle$ of a torsion-free
  word-hyperbolic group $G$, the algorithm computes a number $r$ and
  finite presentations for groups $G_1,\dots, G_k$ (with $k\ge 0$)
  such that
\[
G\cong G_1\ast \dots \ast G_k \ast F_r
\]
and such that each $G_i$ is one-ended.
\end{prop}

\begin{proof}
  A result of Gerasimov~\cite{Ger} shows that there is an algorithm
  that, given a finite presentation of a torsion-free word-hyperbolic
  group $G$, decides whether the group has $1$, $2$ or infinitely many
  ends. If $G$ is two-ended, then $G$ is infinite cyclic. If $G$ is
  one-ended, then it is freely indecomposable.

  If it turns out that $G$ has infinitely many ends we enumerate all
  finite presentations of $G$. This can be done by enumerating all
  possible finite group presentations $W$, all pairs of maps from the
  generators of $G$ to $W$ and from the generating set of $W$ to $G$,
  enumerating the normal closures of the relators of $G$ and $W$ and
  checking if the maps between $G$ and $W$ are group homomorphisms and
  if their compositions define identity maps for $G$ and $W$.
  
  As $G$ is a proper free product we eventually find a presentation of
  type $\langle S_1\cup S_2\,|\,R_1\cup R_2\rangle$ such that any
  $r\in R_i$ is a word in $S_i^{\pm 1}$ for $i=1,2$ and that
  $H=\langle S_1\rangle\neq 1$ and $K=\langle S_2\rangle\neq 1$. As
  both of these properties can be checked we eventually find the
  groups $H$ and $K$ with their finite presentations.
  
  We then iterate the entire process (including determining the number
  of ends) by applying it to each of $H,K$ separately. Eventually we
  will get a collection of infinite cyclic and one-ended groups that
  gives a required decomposition of $G$.
\end{proof}

\section{Virtual fibers and their recognition}

We need a precise description of non-quasiconvex finitely generated
subgroups of groups from $\mathcal M_1$.

\begin{defn}
  Two subgroups $H,K\le G$ of a group $G$, are said to be
  \emph{commensurable} in $G$ if $H\cap K$ has finite index in both
  $H$ and $K$.
  
  Let $G$ be a group and $H\le G$ be a subgroup. The
  \emph{commensurator $\Comm_G(H)$ of $H$ in $G$} is the set of all
  $g\in G$ such that $H$ and $g^{-1}Hg$ are commensurable in $G$.
\end{defn}
The commensurator $\Comm_G(H)$ is easily seen to be a subgroup of $G$
containing~$H$.

\begin{defn}
  Let $G$ be a group from $\mathcal M_1$.
  
  We will say that a subgroup $H\le G$ is a \emph{fiber group} for $G$
  if $H$ is the fundamental group of a closed hyperbolic surface, $H$
  is normal in $G$ and the quotient $G/H\cong \mathbb Z$ is infinite
  cyclic.
  
  We will say that $H\le G$ is a \emph{virtual fiber group} for $G$ if
  $H$ is commensurable in $G$ with a subgroup $H_1$ where $H_1$ is a
  fiber group for some subgroup $G_1\le G$ of finite index in $G$.
\end{defn}

It is obvious that virtual fiber groups are not quasiconvex and, as
it turns out, the converse is also true. The following statement is a
corollary of the work of Canary~\cite{C} on geometrically infinite
ends and of Marden's ``tameness conjecture'' for Kleinian groups
recently proved by Agol~\cite{A} and Calegari--Gabai~\cite{CG}.

\begin{thm}
  Let $G$ be a group from $\mathcal M_1$ and let $H\le G$ be a
  finitely generated subgroup. Then $H$ is quasiconvex in $G$ if and
  only if $H$ is not a virtual fiber group for $G$.
\end{thm}

\begin{rem}\label{rem:vf}
  Note that virtual fibers are locally quasiconvex.  Moreover, it
  follows from the results of Tukia~\cite{Tu} and Gabai~\cite{Ga} that
  if $P$ is virtually a closed hyperbolic surface group and $P$ is
  torsion-free then $P$ is a closed hyperbolic group itself. Thus
  virtual fibers are in fact surface groups. Hence the various
  algorithms discussed in Proposition~\ref{summ}, such as solving the
  uniform membership problem, computing quasiconvexity constants, etc,
  may be considerably improved and sped up in this case (see, for
  example \cite{Sch}).  Moreover, the ``tameness conjecture'' implies
  that the subgroup $H_1$ in the definition of a virtual fiber $H$
  above can be taken as an actual topological fiber group of a finite
  cover $M_1$ fibering over a circle of a closed hyperbolic 3--manifold
  $M$ with the fundamental group $G$.
\end{rem}

\begin{conv}[Computing a subgroup]
  From now on to \emph{compute} a word-hyperbolic subgroup $H$ of a
  word-hyperbolic group $G=\langle S|R\rangle$ shall mean to do all of
  the following:
\begin{enumerate}
\item[(a)] to find a finite generating set for $H$, a finite
  presentation for $H$ on that generating set, a hyperbolicity
  constant for that presentation and the order of $H$;
\item[(b)] if $H$ is known to be quasiconvex in $G$, to find a
  quasiconvexity constant for $H$ in the Cayley graph of the ambient
  group $G$, the index of $H$ in $G$ and a finite state automaton
  accepting the language $L_H(G,S)$ of all short-lex representatives
  over $S$ of elements of $H$;
\item[(c)] if $H$ is known to have finite index in $G$, to find a
  right transveral for $H$ in $G$.
\end{enumerate}
\end{conv}

\begin{prop}\label{fibers}
  There is an algorithm with the following properties.
  
  Suppose we are given a finite presentation $G=\langle S| R\rangle$
  of a group from $\mathcal M_1$, a finite subset $A\subseteq G$ and
  an element $g\in G$.

  Then the algorithm decides whether the subgroup $H:=\langle
  A\rangle\le G$ is quasiconvex in $G$ or not.  If $H$ is quasiconvex,
  the algorithm also produces a quasiconvexity constant for $H$ in
  $X=X(G,S)$, a finite presentation for $H$ on $S$ and the index of
  $H$ in $G$.

  If $H$ is not quasiconvex (and therefore is a virtual fiber group
  for $G$), the algorithm computes a subgroup $G_1$ of finite index in
  $G$ such that $H$ is commensurable with a fiber $H_1$ of $G_1$. The
  algorithm also computes a presentation of $G_1$ as a cyclic
  extension of $H_1$ and it computes a subgroup $H'\le H\cap H_1$ such
  that $H'$ is of finite index in both $H$ and $H_1$.
  
  The algorithm then decides whether or not $g\in H$.
\end{prop}

\begin{proof}

  We will run the ``detection of quasiconvexity'' algorithm from
  Proposition~\ref{summ} in parallel with the following procedure.

  Start enumerating subgroups of finite index in $G$ and computing
  their presentations. For each such subgroup $G_1\le G$ start
  checking if $G_1$ splits as a surface-by-cyclic group
  $G_1=H_1\ltimes \mathbb Z$, where $H_1$ is a surface group. This can
  be done by a ``general enumeration'' algorithm that lists all
  surface-by-cyclic presentations. For each such presentation $W$
  start enumerating maps from the generating set of $W$ to $G_1$ and
  from the generating set of $G_1$ to $W$. For each pair of such maps,
  via enumeration of all the relations in $W$ and $G_1$, start
  checking if they define group homomorphisms and if their
  compositions define identity maps of $W$ and $G_1$ respectively. (We
  refer to this as the ``general enumeration presentation comparison''
  argument.)
  
  If so, then indeed $G_1=H_1\ltimes \mathbb Z$. Note that there is
  an algorithm which, for each $H_1$ above, solves the membership
  problem for $H_1$ in $G$.  Now for each such $H_1$, using the
  algorithm from Lemma~\ref{enum:fi}, start enumerating subgroups
  $H_2$ of $G$ containing $H_1$ as a subgroup of finite index and
  computing finite presentations for such subgroups $H_2$.

  For each of $H_2$ then start checking if $A\subseteq H_2$ (that is
  $H\le H_2$). If yes, check if $H=\langle A\rangle\le H_2$ has finite
  index in the surface group $H_2$. If yes, then $H$ is a virtual
  fiber group for $G$.

  Since $H$ is either quasiconvex or is a virtual fiber group in $G$,
  eventually we will either detect the quasiconvexity of $G$ or will
  discover the fact that $H$ is a virtual fiber for $G$.
  
  In either case it is easy to decide whether or not $g\in H$.
\end{proof}

Thus we see, in particular, that there is a uniform algorithm for
solving the uniform subgroup membership problem for groups from class
$\mathcal M_1$.

\section{The relative rank with respect to virtual fibers}

\begin{defn}[Relative rank]
  Let $G$ be a finitely generated group and let $Z\subseteq G$ be a
  subset. Define the \emph{relative rank} $\rk_Z(G)$ of $G$ with
  respect to $Z$ as the smallest cardinality of a subset $S\subseteq
  G$ such that $G=\langle S\cup Z\rangle$.
\end{defn}

The following statement is an elementary exercise.
\begin{lem}\label{dihed}
  Let $D=\langle a,b| a^2=b^2=1\rangle$ be an infinite dihedral group.

\begin{enumerate}
\item A $k$--tuple $(d_1,\dots, d_k)$ generates $D$ if and only if $k\ge 2$
  and this $k$--tuple is Nielsen-equivalent to $(a,b,1,\dots, 1)$.
\item Suppose $s$ has order two in $D$ and $d_1,\dots, d_k\in \langle
  ab\rangle$. Then $(s,d_1,\dots, d_k)$ generates $D$ if and only
  $(d_1,\dots, d_k)$ generates $\langle ab\rangle$, that is, if and
  only if, $(d_1,\dots, d_k)$ is Nielsen-equivalent to $(ab, 1,\dots,
  1)$ in $D$.
\end{enumerate}
\end{lem}

\begin{prop}\label{vc}
  There is an algorithm that, given a finite presentation of a
  virtually cyclic group $G$ and a finite subset $Z\subseteq G$, finds
  $\rk_Z(G)$.
\end{prop}

\begin{proof}
  We first compute the order of $G$ and the order of $L:=\langle
  Z\rangle$. If $G$ is finite or $L$ is infinite (that is
  $[G:L]<\infty$), the problem easily reduces to the relative rank
  problem for finite groups.
  
  Therefore we may assume that $G$ is infinite and $L$ is finite.
  Then, by a well-known fact about virtually cyclic groups, $G$
  possesses a finite normal subgroup $N$ such that $\overline{G}:=G/N$
  is either infinite cyclic or infinite dihedral. By a general
  enumeration argument we can find such a subgroup $N$ and determine
  which of these two possibilities occurs.

{\bf Case 1}\qua The group $\overline G$ is infinite cyclic.

In this case clearly $N$ is the set of all torsion elements of $G$ and
hence $L\le N$.  Find an element $x\in G$ such that $\overline
G=\langle \overline x\rangle$.

By performing the Euclidean algorithm modulo $N$ we see that if a tuple
$(g_1,\dots, g_k,Z)$ generates $G$ then the $k$--tuple $(g_1,\dots, g_k)$
is Nielsen-equivalent to a $k$--tuple of the form $(xn_1,n_2,\dots, n_k)$
where $n_i\in N$.

To decide if $\rk_{Z}(G)\le k$ we enumerate all tuples
$(xn_1,n_2,\dots, n_k,Z)$, where $n_i\in N$, and check if at least one
of them generates $G$.

{\bf Case 2}\qua The group $\overline G$ is infinite dihedral.

Then $\overline G=\langle \overline s, \overline t~|~\overline{s}^2=\overline{t}^2=1\rangle$. We compute $s,t\in G$ that
map to $\overline s,\overline t$ accordingly. The infinite cyclic
subgroup $\langle \overline{st}\rangle$ has index two in $\overline
G$. Let $H\le G$ be the full preimage of $\langle
\overline{st}\rangle$ in $G$, so that $H=\langle st, N\rangle$ has
index two in $G$.

Note that $\overline L$ is either trivial or cyclic of order two.
Suppose first that $\overline L$ is trivial, that is $Z\subseteq N$.
Then by Lemma~\ref{dihed} if $(g_1,\dots, g_k,Z)$ generates $G$ then
$k\ge 2$ and $(g_1,\dots, g_k)$ is Nielsen-equivalent to a $k$--tuple of
the form $(sn_1,tn_2,n_3,\dots, n_k)$ where $n_i\in N$. To decide if
$\rk_{Z}(G)\le k$ we enumerate all tuples $(sn_1,tn_2,n_3\dots,
n_k,Z)$, where $n_i\in N$, and check if at least one of them generates
$G$.

Suppose now that $\overline L$ is cyclic of order two, that is there
is some element $z\in Z$ such that $\overline z$ has order two. In
this case $L\cap H\le N$. Also, if $z'\in Z$ is different from $z$,
then either $z'\in H$ (and hence $z'\in N$) or
$\overline{z'}=\overline z$.

Suppose $(g_1,\dots, g_k,Z)$ generates $G$.  After multiplying $g_i$
by $z^{-1}$ if necessary, we obtain a tuple $(g_1',\dots, g_k',Z)$
generating $G$ such that all $g_i'\in H$. Note that $(\overline
{g_1'},\dots, \overline{g_k'},\overline z)$ generates $\overline G$
since for each $z'\in Z$ either $\overline{z'}=\overline z$ or
$\overline{z'}=1$.  Now Lemma~\ref{dihed} implies that $(g_1',\dots,
g_k')$ is Nielsen-equivalent to $(stn_1,n_2,\dots, n_k)$ where $n_i\in
N$.

To decide if $\rk_{Z}(G)\le k$ we enumerate all tuples
$(stn_1,n_2,\dots, n_k,Z)$ where $n_i\in N$, and check if at least one
of them generates $G$.
\end{proof}

\begin{rem}\label{rem:vc}
  The proof of Proposition~\ref{vc} actually shows that if $G$ is
  virtually cyclic and $Z$ is a tuple generating a finite subgroup of
  $G$ then for every $k\ge 1$ there is a finite set $E$ of tuples
  $(f_1,\dots, f_k,Z)$ such that every tuple $(g_1,\dots, g_k,Z)$
  generating $G$ is Nielsen-equivalent to a tuple from $E$.
\end{rem}

\begin{lem}\label{conjug}
  Let $G\in \mathcal M_1$, let $G_1$ be a subgroup of finite index in
  $G$ and let $H\le G_1$ be a fiber subgroup for $G_1$. Thus
  $G_1=H\ltimes \langle t\rangle$. Let $P\le H$ be a subgroup of
  finite index in $H$ that is normal in $G_1$.
  
  Then the following hold:
\begin{enumerate}
\item For each subgroup $K$ of $G$ such that $P\le K$, either $K$ has
  finite index in $G$ or $K$ contains $P$ as a subgroup of finite
  index.
\item The set of subgroups conjugate to $P$ in $G$ is finite.
\item For each $g\in G$ either $P$ has finite index in $L:=\langle
  P,g^{-1}P g\rangle$ (in which case $g\in \Comm_G(P)$) or $L$ has
  finite index in $G$.
\item Either there is $g\in G-G_1$ such that $\langle
  P,g^{-1}Pg\rangle$ has finite index in $G$ or for every $g\in G$ the
  subgroup $P$ has finite index in $\langle P,g^{-1}Pg\rangle$ (in
  which case $G=\Comm_G(P)$).
\end{enumerate}
\end{lem}

\begin{proof}
  Let $g_1=1, \dots, g_n$ be a right transversal for $G_1$ in $G$.
  Then every element of $G$ is uniquely expressible in the form
  $g=ht^jg_i$ where $1\le i\le n$, $j\in \mathbb Z$ and $h\in H$.
  Recall that $P$ is normal in $G_1$.  Hence
  $g^{-1}Pg=g_i^{-1}t^{-j}h^{-1}Pht^jg_i=g_i^{-1}Pg_i$. This shows
  that the set of conjugates of $P$ in $G$ is finite and part (2) of
  the lemma is established.
  
  Let $K\le G$ be a subgroup such that $P\le K$ and suppose that $P$
  has infinite index in $K$. Let $K_1=K\cap G_1$. Then $P\le K_1\le K$
  and $K_1$ has finite index in $K$. Therefore $P$ has infinite index
  in $K_1$. Since $K_1\le G_1$, every element of $K_1$ can be written
  in the form $t^jh$ where $j\in \mathbb Z$ and $h\in H$. Assume first
  that there is some element $a\in K_1$ of the form $a=t^jh$, where
  $j\ne 0$. Then it is easy to see that $\langle P, a\rangle$ has
  finite index in $G_1$ and hence in $G$, as required.
  
  Suppose now that there is no element $a\in K_1$ of the form $a=t^jh$
  with $j\ne 0$. This implies that $P\le K_1\le H$, contradicting our
  assumption that $[K_1:P]=\infty$. This verifies part (1) of the
  lemma.

  Part (3) follows immediately from part (1).  In turn, part (3)
  immediately implies part (4).
\end{proof}

\begin{lem}\label{cases}
  There is an algorithm with the following properties. Suppose we are
  given a group $G\in \mathcal M_1$, a subgroup of finite index
  $G_1\le G$ and a fiber group $H$ of $G_1$ so that $G_1=H\ltimes
  \langle t\rangle$. Suppose we are also given a subgroup $P\le H$ of
  finite index in $H$ that is normal in $G_1$.

  The algorithm decides if there is $g\in G-G_1$ such that $\langle
  P,g^{-1}Pg\rangle$ has finite index in $G$ or if for every $g\in G$
  the subgroup $P$ has finite index in $\langle P,g^{-1}Pg\rangle$.
  
  In the former case the algorithm computes all the (finitely many)
  subgroups $\langle P,g^{-1}Pg\rangle$, where $g\in G$, that have
  finite index in $G$. It then computes their intersection and finds a
  normal subgroup $N\le G$ of finite index in $G$ that is contained in
  that intersection and such that $N\le G_1$.
  
  In the latter case the algorithm computes the intersection $N_1$ of
  all conjugates of $P$ in $G$. Thus $N_1$ is a normal subgroup of $G$
  that has finite index in $P$ and in $H$.
\end{lem}

\begin{proof}
  We first compute a right transversal $g_1=1, g_2, \dots g_n$ for
  $G_1$ in $G$. We have already seen in the proof of
  Lemma~\ref{conjug} that each conjugate of $P$ in $G$ has the form
  $g_i^{-1}Pg_i$. Put $L_i:=\langle P, g_i^{-1}Pg_i\rangle\le G$. For
  each $i$ either $P$ has finite index in $L_i$ or $L_i$ has finite
  index in $G$. For each $i$ we run the algorithm provided by
  Lemma~\ref{enum:fi} in parallel with the Todd--Coxeter coset
  enumeration algorithm for $L_i$ to decide which of these
  alternatives holds.
  
  In particular, we find all of the $L_i$ (if any) that have finite
  index in $G$. If there is at least one such $L_i$, we find their
  intersection $L$, and the intersection $L'=L\cap G_1$. Then $L'$
  still has finite index in $G$. We then find the intersection of all
  conjugates of $L'$ in $G$ and denote it by $N$. Clearly $N$ is
  normal of finite index in $G$, also $N$ is a subgroup of $G_1$ and
  $N$ is contained in all those $L_i$ that have finite index in $G$,
  as required.

  Suppose now that all $L_i$ turn out to contain $P$ as a subgroup of
  finite index. Thus each $L_i$ is a virtual fiber for $G$. For each
  of $L_i$ we compute its finite presentation. Recall that all
  conjugates of $P$ in $G$ are of the form $P_i:=g_i^{-1}Pg_i$. For
  each $i$, operating inside $L_i$ we compute a generating set for
  $P_i':=P_i\cap P$ and rewrite it as the set of words in the
  generators of $P$. This is possible by Proposition~\ref{summ} since
  $L_i$ is a surface group and thus is locally quasiconvex. Finally,
  operating inside the surface group $H$, we compute the subgroup
$$N_1:=\cap_{i=1}^n P_i'=\cap_{i=1}^n P_i'=\underset{g\in G}{\cap} g^{-1}Pg.$$
Then $N_1$ is a subgroup of finite index in $H$ and $N_1$ is normal in
$G$, as required.
\end{proof}

\begin{prop}\label{rel}
  There is an algorithm with the following properties.  Suppose we are
  given a group $G\in \mathcal M_1$, a subgroup of finite index
  $G_1\le G$, a fiber group $H$ of $G_1$ so that $G_1=H\ltimes \langle
  t\rangle$ and a virtual fiber group $H_1$ commensurable with $H$.
  
  The algorithm computes the relative rank $\rk_{H_1}(G)$ of $G$ with
  respect to $H_1$.
\end{prop}

\begin{proof}
  
  First we compute a subgroup $P\le H$ of finite index in $H$ such
  that $P\le H_1$ and such that $P$ is normal in $G_1$. Such a
  subgroup obviously exists since $H\cap H_1$ has finite index in $H$
  and hence contains a subgroup $P'$ such that $P'$ is a
  characteristic subgroup of finite index in $H$. Since conjugation by
  $t$ induces an automorphism of $H$, such a subgroup $P'$ will be
  normal in $G_1$.
  
  We can algorithmically find some subgroup $P$ with the required
  properties as follows. First, compute the automorphism $\phi$ of $H$
  induced by conjugation by $t$. That is, for each generator $x$ of
  $H$ express $t^{-1}xt$ as a word in the generators of $H$. Then
  start enumerating subgroups of finite index $P$ in $H$. For each
  such $P$ start checking if $P$ is normal in $H$, if $P$ is contained
  in $H_1$ and if $\phi(P)=P$. If yes, then $P$ is as required, and we
  terminate the process. As we observed above, this algorithm will
  necessarily terminate since some $P$ with required properties does
  exist.

  We then use the algorithm from Lemma~\ref{cases} and decide if there
  is $g\in G-G_1$ such that $\langle P,g^{-1}Pg\rangle$ has finite
  index in $G$ or if for every $g\in G$ the subgroup $P$ has finite
  index in $\langle P,g^{-1}Pg\rangle$.
  
  Suppose first that the former occurs. Then we compute the normal
  subgroup $N$ of $G$ of finite index in $G$ defined in
  Lemma~\ref{cases}. Denote $\overline G=G/N$.
  
  {\bf Claim}\qua For any set $B\subseteq G$ we have $\langle
  H_1,B\rangle=G$ if and only if $\langle \overline H_1, \overline
  B\rangle=\overline G$. Indeed, recall that $N\le G_1$ is normal of
  finite index in $G$. Since $P$ is normal in $G_1$, for every $g\in
  G$ and every $n\in N$ we have $g^{-1}n^{-1}Png=g^{-1}Pg$ and hence
  $\langle P,g^{-1}Pg\rangle=\langle P,g^{-1}n^{-1}Png\rangle$.
  
  Suppose $\langle \overline H_1, \overline B\rangle=\overline G$.
  Choose $g\in G$ such that $\langle P,g^{-1}Pg\rangle$ has finite
  index in $G$. Then there is $n\in N$ such that $gn\in \langle
  H_1,B\rangle$. Hence $\langle P,g^{-1}Pg\rangle=\langle
  P,g^{-1}n^{-1}Png\rangle$ has finite index in $G$ and therefore
  contains $N$.  But both $P$ and $gn$ are contained in $\langle H_1,
  B\rangle$.  Hence $N\le\langle H_1, B\rangle$. Since $\langle
  \overline H_1, \overline B\rangle=\overline G$, this implies that
  $\langle H_1,B\rangle=G$. This verifies the Claim.
  
  We now compute the image $\overline H_1$ of $H_1$ in the finite
  group $\overline G$ and solve the relative rank problem for
  $\overline G$ with respect to $\overline H_1$.  By the Claim
  $\rk_{H_1}(G)=\rk_{\overline H_1} (\overline G)$.

  Suppose now that the second alternative of Lemma~\ref{cases} occurs.
  Then, as in Lemma~\ref{cases} we compute the intersection $N_1$ of
  all conjugates of $P$ in $G$. This is a normal subgroup of $G$ which
  has finite index in $H$.
  
  Then the group $\widehat G:=G/N_1$ is virtually cyclic. Since $N_1\le
  H_1$, it is clear that $\rk_{H_1}(G)=\rk_{\widehat H_1}(\widehat G)$ where
  $\widehat H_1$ is the image of $H_1$ in $\widehat G$. We use then use the
  algorithm provided by Proposition~\ref{vc} to compute $\rk_{\widehat
  H_1}(\widehat G)$.
\end{proof}

\section{Main technical tool}

In this section we will state the main technical tool needed to prove
our main result. This tool was obtained by Kapovich and Weidmann
in~\cite{KW04}.

If $n\ge 0$ is an integer, for an $n$--tuple $T=(g_1,\dots, g_n)$ of
elements of a group $G$ denote $n=l(T)$ and call $n$ the \emph{length
  of $T$}.

Let us recall the notion of Nielsen equivalence:

\begin{defn}[Nielsen equivalence]
  Let $T=(g_1,\dots ,g_n)\in G^n$ be an $n$--tuple of elements of a group $G$.
  The following moves are called \emph{elementary Nielsen moves} on
  $T$:

\begin{enumerate}
\item[(N1)] For some $i, 1\le i\le n$ replace $g_i$ by $g_i^{-1}$ in
  $T$.
\item[(N2)] For some $i\ne j$, $1\le i,j\le n$ replace $g_i$ by
  $g_ig_j$ in $T$.
\item[(N3)] For some $i\ne j$, $1\le i,j\le n$ interchange $g_i$ and
  $g_j$ in $T$.
\end{enumerate}

We say that $T=(g_1,\dots, g_n)\in G^n$ and $T'=(f_1,\dots , f_n)\in
G^n$ are \emph{Nielsen-equivalent in $G$}, if there is a chain of
elementary Nielsen moves which transforms $T$ to $T'$.
\end{defn}

\begin{defn}[Partitioned tuple]
  Let $G$ be a group. A \emph{partitioned tuple} for $G$ is a tuple
  $M=(Y_1,\dots,Y_s;T)$ where $s\ge 0$, and where

\begin{enumerate}
\item[(a)] each of $Y_i,T$ is a tuple of elements of $G$;
\item[(b)] either $s>0$ or $l(T)>0$;
\item[(c)] we have $\langle Y_i\rangle \ne 1$ for each $i>0$.
\end{enumerate}

We call $l(Y_1)+\dots +l(Y_s)+l(T)$ the \emph{length} of $M$ and
denote it by $l(M)$.  We call the $l(M)$--tuple of elements of $G$
obtained by concatenating the tuples $Y_1,\dots, Y_s,T$ the
\emph{underlying tuple of $M$}.
\end{defn}

Thus $(;T)$ (where $l(T)>0$) and $(Y_1; )$ (where $\langle
Y_1\rangle\ne 1$) are examples of partitioned tuples.

\begin{defn}[Elementary Moves]
  Let $M=(Y_1,\dots,Y_s;T)$ be a partitio\-ned tuple for a group $G$.
  Let $T=(t_1,\dots, t_m)$.  The \emph{elementary moves} on $M$ are
  the following:
\begin{enumerate}
\item Replace some $Y_i$ by $g^{-1}Y_ig$ where $g\in \langle
  (\cup_{j\ne i} Y_j)\cup T\rangle$.
\item Replace some entry $t_i$ in $T$ by $ut_iu'$ where \[u,u'\in
  \langle (\cup_{j=1}^s Y_j)\cup\{t_1, \dots, t_{i-1}, t_{i+1},\dots,
  t_m \} \rangle.\]
\end{enumerate}
\end{defn}

\begin{defn}[Equivalence of partitioned tuples]
  Two partitioned tuples\linebreak $M=(Y_1,\dots,Y_s;T)$ and
  $M'=(Y_1',\dots,Y_s';T')$ for a group $G$ are \emph{equivalent} if
  there exists a chain of elementary moves taking $M$ to $M'$.
\end{defn}
It is easy to see that the above definition gives an equivalence
relation on the set of partitioned tuples. Moreover, if $M$ and $M'$
are equivalent partitioned tuples then $l(M)=l(M')$ and the underlying
tuples of $M, M'$ are Nielsen-equivalent in~$G$.

\begin{defn}[Complexity of a partitioned tuple]
  We define the complexity of a partitioned tuple $M=(Y_1,\dots
  ,Y_n;T)$ with $T=(t_1,\ldots ,t_m)$ to be the pair $(m,n)\in\mathbb
  N^2$. We define an order on $\mathbb N^2$ by setting $(m,n)\le
  (m',n')$ if $m<m$ or if $m=m'$ and $n\le n'$. This clearly gives a
  well-ordering on $\mathbb N^2$.
\end{defn}

\begin{notation}[Invariant sets]
  Let $G$ be a nonelementary torsion-free word-hyperbolic group $G$
  with a finite generating set $S$. Let $X$ be the Cayley graph of $G$
  with respect to $S$. Let $\delta\ge 1$ be an integer such that $X$
  is $\delta$--hyperbolic.
  
  Let $U\le G$ be a nontrivial subgroup (which is therefore infinite).
  Let $\Lambda U\subseteq \partial X$ be the limit set of $U$ in $X$.
  
  Let $E(U)$ be the set of all $x\in X$ such that for some $g\in U,
  g\ne 1$ we have $d(x,gx)\le 100\delta$. Let $Z(U)$ be the weak
  convex hull of $E(U)\cup \Lambda U$, that is $Z(U)$ is the union of
  all geodesics in $X$ with both endpoints in $E(U)\cup \Lambda U$.
  Finally, let $X(U)$ denote the closure of $Z(U)$ in $X$.
\end{notation}

Note that the definitions of the above $U$--invariant subsets of $X$
depend on the choice of $\delta$.

The following is a corollary of \cite[Theorem~2.4]{KW04}. This is
the main technical tool required for our proofs in this paper.

\begin{thm}\label{tool}
  For every integer $k\ge 1$ there exists a computable constant
  $L=L(k)\ge 0$ with the following property.
  
  Let $G$ be a nonelementary torsion-free word-hyperbolic group with a
  finite generating set $S$. Let $X$ be the Cayley graph of $G$ with
  respect to $S$. Let $\delta\ge 1$ be an integer such that $X$ is
  $\delta$--hyperbolic.
  
  Let $M=(Y_1,\dots, Y_s; T)$ by a partitioned tuple for $G$ with
  $l(M)=k$.  Let $U$ be the subgroup generated by the underlying tuple
  of $M$.  Then either
\[
U=\langle Y_1\rangle \ast \dots \langle Y_s\rangle \ast F(T)
\]
or $M$ is equivalent to a partitioned tuple $M'=(Y_1',\dots, Y_s';
T')$ such that, denoting $U_i'=\langle Y_i'\rangle$, one of the
following occurs:

\begin{enumerate}
  
\item There are some $i\ne j$ such that $d(X(U_i'), X(U_j'))\le \delta
  L(k)$.
  
\item There is some $i, 1\le i\le s$ and some element $g$ of $T'$ such
  that \[ d(X(U_i'), gX(U_i'))\le \delta L(k).\]
  
\item There is some element $g$ of $T'$ such that $g$ is conjugate in
  $G$ to an element of length at most $\delta L(k)$.

\end{enumerate}

\end{thm}

\section{Generator transfer process}

If $G$ is a group with a finite generating set $S$ and if $g\in G$, we
will denote by $|g|_S$ the $S$--geodesic length of $g$.  Thus
$|g|_S=d(g,1)$ in $X(G,S)$.  All the constants in this section are
assumed to be monotone non-decreasing in each of their integer
arguments.

We recall here some results of Kapovich--Weidmann~\cite{KW04}.

\begin{lem}\label{close} There
  exists a computable constant $c_1=c_1(G,S,\delta,\epsilon)>0$ such
  that the following holds:
  
  Suppose $G=\langle S|R\rangle$ is a finite presentation of a
  torsion-free word-hyperbolic group $G$ and that $\delta\ge 0$ is an
  integer such that the Cayley graph $X(G,S)$ is $\delta$--hyperbolic.
  Let $\epsilon\ge 0$ be an integer and let $U\le G$ be an infinite
  subgroup such that $U$ is $\epsilon$--quasiconvex in $X=X(G,S)$. Then
  $U$ and $X(U)$ are $c_1$--Hausdorff close (where $X(U)$ is defined
  relative $\delta$).
\end{lem}
Lemma~\ref{close} follows from \cite[Lemma~10.3 and Remark~10.9]{KW04}.

\begin{lem}\label{hyp4} For any integers
  $K\ge 0$ and $n\ge 1$ there is a computable constant
  $c_2=c_2(G,S,\delta,n,K)$ with the following property.  Suppose
  $G=\langle S|R\rangle$ is a finite presentation of a group $G$ from
  the class $\mathcal M_1$. Suppose $\delta\ge 1$ is an integer such
  that the Cayley graph $X=X(G,S)$ is $\delta$--hyperbolic. Suppose
  $U\le G$ is a non-trivial quasiconvex subgroup $U\le G$ generated by
  a set $Y$ with at most $n$ elements such that $Y$ is contained in
  the $K$--ball around $1$ in $\Gamma(G,S)$. Then the sets $U$ and
  $X(U)$ are $c_2$--Hausdorff close (where $X(U)$ is defined relative
  to $\delta$).
\end{lem}

\begin{proof}
  Let $E$ be the maximum of quasiconvexity constants of infinite
  quasiconvex subgroups generated by subsets with at most $n$--elements
  from the $K$--ball around $1$ in $\Gamma(G,S)$.  Put
  $c_2:=c_1(G,S,\delta,E)$. Then $c_2$ clearly satisfies the
  requirements of the lemma. In order to see that $c_2$ is computable
  it suffices to see that $E$ is computable.  For each $n$--tuple $Y$
  of words of length at most $K$ over $S$ apply the algorithm of
  Proposition~\ref{fibers} to decide whether or not $H=\langle
  Y\rangle$ is quasiconvex in $G$, and if yes, to compute a
  quasiconvexity constant for $H$. Then $E$ is the maximum of all the
  quasiconvexity constants obtained in this way and hence $E$ is
  computable.
\end{proof}

\begin{prop}\label{hyp5}
  There exist a computable constant $c=c(G,S,\delta, n_1,n_2,K)$ with
  the following properties.
  
  Suppose $G=\langle S|R\rangle$ is a finite presentation of a group
  $G$ from the class $\mathcal M_1$. Suppose $\delta\ge 1$ is an
  integer such that the Cayley graph $X=X(G,S)$ is
  $\delta$--hyperbolic. Suppose that $K,n_1,n_2\ge 1$ are integers.

  Suppose $U_1=\langle Y_1\rangle$ and $U_2=\langle gY_2g^{-1}\rangle$
  are two quasiconvex subgroups of $G$, where $g\in G$,
  $Y_1=(y_1,\ldots ,y_{n_1})\in G^{n_1}$, $Y_2=(y_1',\ldots
  ,y_{n_2}')\in G^{n_2}$, $|y_i|_S\le K$ for $1\le i\le n_1$ and
  $|y_i'|_S\le K$ for $1\le i\le n_2$. Suppose also that
  $d(X(U_1),X(U_2))\le K$.
  
  Then the $(n_1{+}n_2)$--tuple $(y_1,\ldots ,y_{n_1},gy_1'g^{-1},\ldots
  ,gy_{n_2}'g^{-1})$ is Nielsen-equival\-ent to a tuple conjugate in $G$
  to $(y_1,\ldots ,y_{n_1+n_2})$ where $|y_i|_S\le c$ for $1\le i\le
  n_1+n_2$.
\end{prop}

\begin{proof}
  It is proved in \cite{KW04} that $c=4c_2(G,S,\delta,n,K)+3K$
  satisfies the requirements of the proposition. The computability of
  $c_4$ now follows from the computability of $c_2$ established in
  Lemma~\ref{hyp4}.
\end{proof}

\begin{prop}\label{hyp6}
  There is a computable constant $c'=c'(G,S,\delta, n_1,n_2,K)$
  with the following properties.
  
  Suppose $G=\langle S|R\rangle$ is a finite presentation of a group
  $G$ from the class $\mathcal M_1$. Suppose $\delta\ge 1$ is an
  integer such that the Cayley graph $X=X(G,S)$ is
  $\delta$--hyperbolic.  Suppose that $K,n\ge 1$ are integers.
  
  Let $U=\langle Y\rangle$ where $Y=(y_1,\ldots ,y_{n})\in G^{n}$ and
  $|y_i|\le K$ for $1\le i\le n$. Suppose that $U$ is quasiconvex in
  $G$ and that $d(X(U),gX(U))\le K$.
  
  Then the tuple $(y_1,\ldots ,y_{n},g)$ is Nielsen-equivalent to a
  tuple $(y_1,\ldots ,y_n,y_{n+1})$ such that $|y_i|_S\le c'$ for
  $1\le i\le n+1$.
\end{prop}

\begin{proof}
  It is proved in \cite{KW04} that $c'=2c_2(G,S,\delta,n,K)+3K$
  satisfies the requirements of the proposition.  By Lemma~\ref{hyp4}
  $c_2$ is computable and therefore $c'$ is also computable.
\end{proof}

\begin{thm}\label{constC}
  There exists a computable integer constant $C=C(G,S,\delta,k)>0$
  with the following properties.
  
  Suppose $G=\langle S|R\rangle$ be a finite presentation of a group
  $G$ from the class $\mathcal M_1$. Suppose $\delta\ge 1$ is an
  integer such that the Cayley graph $X=X(G,S)$ of $G$ with respect to
  $S$ is $\delta$--hyperbolic. Suppose $k\ge 1$ is an integer.

  Then for any $k$--tuple $T=(g_1,g_2,\dots, g_k)$ generating $G$ one
  of the following holds:

\begin{enumerate}
\item There exists a $k$--tuple $T'=(g_1',g_2',\dots, g_k')$
  Nielsen-equivalent to $T$ such that $|g_i'|_S\le C$ for $i=1,\dots,
  k$.
\item There exists a $k$--tuple $T'=(g_1',g_2',\dots, g_k')$
  Nielsen-equivalent to $T$ such that for some $j<k$ the set
  $(g_1',\dots, g_j')$ generates a virtual fiber of $G$ and
  $|g_i'|_S\le C$ for $i=1,\dots, j$.
\end{enumerate}
\end{thm}

\begin{proof}
  Let $G=\langle S|R\rangle$ be a finite presentation of a group from
  $\mathcal M_1$ and let $\delta\ge 1$ be an integer such that
  $X=X(G,S)$ is $\delta$--hyperbolic.
  
  We will prove the result by induction on $k$. For $k=1$ there is
  nothing to prove, since $G$ is not cyclic. Suppose now that $k\ge 2$
  and that the statement has been proved for all $1\le k'<k$.
  
  We define the constants $R(i)=R(k,i)$ for $i=1,\dots, k$ inductively
  as follows.  Put $R(1):=\delta L(k)$, where $L(k)$ is the constant
  provided by Theorem~\ref{tool}.  For $1<i\le k$ put
$$R(i)=\max\Big\{R(i{-}1),\,c'\left(G,S,\delta,i{-}1,R(i{-}1)\right),\,
  \max_{\stackrel{\scriptstyle p+j=i}{\scriptstyle p,j\ge 1}}
  c\left(G,S,\delta,p,j,R(i{-}1)\right)\Big\}.$$
Note that the constants $R(k,i)$ are algorithmically computable.

Let $M=(Y_1,\dots, Y_s;H)$ be a partitioned tuple of length $k$.  Let
$n_i$ be the length of $Y_i$.  We will say that $M$ is \emph{good} if
one of the following holds:

(a)\qua each $Y_j$ generates a nontrivial quasiconvex subgroup of $G$
and each $Y_j$ is conjugate in $G$ to a tuple contained in the ball of
radius $R(n_i)$ in $X$.

(b)\qua there is some $Y_j$ which generates a virtual fiber of $G$ and
such that $Y_i$ is conjugate to a tuple contained in the ball of radius
$R(n_i)$ in $X$.

Put $C(G,S,\delta,k):=\max\{ C(G,S,\delta,k-1), R(k)\}$.
Suppose now that $G$ is generated by a $k$--tuple $T_0=(g_1,\dots,
g_k)$. If $T_0$ is Nielsen-equivalent in $G$ to some tuple containing
$1\in G$, the statement follows from the inductive hypothesis.
Therefore we will assume that every $k$--tuple Nielsen-equivalent to
$T_0$ consists of nontrivial elements of $G$.

Consider a partitioned tuple $M_0:=(; T_0)$. Note that $M_0$ is good
and that the underlying tuple of $M_0$ is $T_0$.  Among all the good
partitioned tuples with the underlying tuple Nielsen-equivalent to
$T_0$ choose a partitioned tuple $M$ of minimal complexity.

{\bf Claim}\qua
We claim that either $M$ has the form $M=(Y_1;)$ or part (b) of the
definition of a good partitioned tuple applies to $M$.

Indeed, suppose not.  Then $M=(Y_1,\dots, Y_s; T)$ where either $s\ge
2$ or $l(T)>0$ and where each $Y_j$ generates a quasiconvex subgroup
of $G$.  Denote $n_i=l(Y_i)$ and $n=l(T)$. Recall that the underlying
tuple of $M$ is Nielsen-equivalent to $T_0$ and hence $n_1+\dots+
n_s+n=k$.

Since $M$ generates a freely indecomposable group $G$, there exists a
partitioned tuple $M'=(Y_1',\dots Y_s';T')$ equivalent to $M$ and such
that one of the three cases of Theorem~\ref{tool} applies to $M'$.
Recall that by definition of equivalence of partitioned tuples each
$Y_i'$ is conjugate to $Y_i$ and the underlying tuple of $M'$ is
Nielsen-equivalent to the underlying tuple of $M$.  In particular $M'$
is good, since $M$ was good. Put $U_i=\langle Y_i'\rangle$.

Suppose that Case (1) of Theorem~\ref{tool} applies to $M'$.  Without
loss of generality we may assume that $d(X(U_1),X(U_2))\le
L(k)\delta$.  Recall that since $M$ and $M'$ are good, each $Y_i'$ is
conjugate to a tuple contained in the ball of radius $R(n_i)$ in $X$.
Then by Proposition~\ref{hyp5} the tuple $(Y_1',Y_2')$ is
Nielsen-equivalent to an $(n_1{+}n_2)$--tuple $Y_0''$ such that $Y_1''$
is conjugate to a set contained in the ball of radius
$c=c(G,S,\delta,n_1,n_2, \max\{L(k)\delta, R(n_1), R(n_2)\})$ in $X$.
Recall that by definition of $R(i)$ we have $R(1)=L(k)\delta\le
R(n_1), R(n_2)\le R(n_1+n_2-1)$. Therefore $Y_1''$ is conjugate to a
set contained in the ball of radius $R(n_1+n_2)$.  Consider now a
partitioned tuple $M'':=(Y_1'',Y_3',\dots, Y_s';T')$. Then $M''$ is
good and the underlying tuple of $M''$ is Nielsen-equivalent to $T_0$.
On the other hand the complexity of $M''$ is smaller than the
complexity of $M$, which contradicts the minimal choice of $M$.

Suppose now that Case (2) of Theorem~\ref{tool} applies to $M'$.
Without loss of generality we may assume that for some entry $g$ of
$T'$ we have $d(X(U_1), gX(U_1))\le L(k)\delta$. Then by
Proposition~\ref{hyp6} the $(n_1{+}1)$--tuple $(Y_1',h)$ is
Nielsen-equivalent to an $(n_1{+}1)$--tuple $Y_1''$ which is conjugate to
a subset contained in the ball of radius
$c'=c'(G,S,\delta,n_1,\max\{R(n_1), L(k)\delta\})\le R(n_1+1)$ in $X$.
Put $M'':=(Y_1'',Y_2',\dots, Y_s'; T'')$ where $T''$ is obtained from
$T'$ by deleting the entry $g$. Then $M''$ is good and the underlying
tuple of $M''$ is Nielsen-equivalent to $T_0$. On the other hand the
complexity of $M''$ is smaller than that of $M$, which contradicts the
minimal choice of $M$.

Suppose that Case (3) of Theorem~\ref{tool} applies to $M'$. Then $T'$
contains an element $g$ that is conjugate in $G$ to an element of
length at most $L(k)\delta$.  Put $M'':=(Y_1',\dots Y_s', Y_{s+1};
T'')$ where $Y_{s+1}=(g)$ and where $T''$ is obtained from $T'$ by
deleting the entry $g$. Then $M''$ is good and the underlying tuples
of $M'$ and $M''$ coincide. In particular, the underlying tuple of
$M''$ is Nielsen-equivalent to $T_0$. On the other hand the complexity
of $M''$ is smaller than that of $M$, which contradicts the minimal
choice of $M$.

Thus we have verified the Claim. This immediately implies the
statement of the theorem.
\end{proof}

\section{The rank problem}

In this section we give a proof of Theorem~\ref{main}:

\begin{thm}\label{m1}
  There is an algorithm that, given a finite presentation $G=\langle
  S|R\rangle$ of a group $G$ from the class $\mathcal M_1$, computes
  the rank of $G$.
\end{thm}

\begin{proof}
  Let $G=\langle S|R\rangle$ be a finite presentation of $G\in
  \mathcal M_1$.  Denote $X=X(G,S)$, the Cayley graph of $G$ with
  respect to $S$. First we compute an integer $\delta\ge 1$ such that
  $X$ is $\delta$--hyperbolic. Let $k=\#S$. We next compute the
  constant $C=C(G,S,\delta,k)>0$ provided by Theorem~\ref{constC}.  We
  then enumerate all $p$--tuples $T$ of words over $S$ with $1\le p\le
  k$, where every word has length at most $C$. Let $\mathcal T$ be the
  set of all such tuples.

  For each tuple $T\in \mathcal T$ we use Proposition~\ref{fibers} to
  decide if $H_T:=\langle T\rangle$ is quasiconvex in $G$ or if $H_T$
  is a virtual fiber of $G$. If $H_T$ is a virtual fiber of $G$ we
  compute the rank of $H_T$ and, using Proposition~\ref{rel}, the
  relative rank $\rk_{H_T}(G)$ of $G$ with respect to $H_T$. If $H_T$
  is quasiconvex in $G$, we check whether or not $H_T=G$, that is
  $\langle T\rangle=G$.
  
  Let $p\le k$. Theorem~\ref{constC} implies that $\rk(G)\le p$ if and
  only if either there exists a $p$--tuple $T\in \mathcal T$ such that
  $G=\langle T\rangle$ or there exists a tuple $T\in \mathcal T$
  generating a virtual fiber $H_T$ such that $\rk(H_T)+\rk_{H_T}(G)\le
  p$.
  
  Therefore for each $0\le p\le k$ we can decide whether or not
  $\rk(G)\le p$. The smallest $0\le p\le k$ such that $\rk(G)\le p$ is
  the rank of $G$.
\end{proof}

\begin{thm}\label{m}
  There is an algorithm that, given a finite presentation $G=\langle
  S|R\rangle$ of a group $G$ from the class $\mathcal M$, computes the
  rank of $G$.
\end{thm}
\begin{proof}
  We first apply the algorithm from Proposition~\ref{ger} and compute
  a decomposition
\[
G=G_1\ast \dots \ast G_s \ast F_r
\]
where each $G_i$ is one-ended and is given to us by an explicit finite
presentation.

By Proposition~\ref{thur} each $G_i\in \mathcal M_1$. We use the
algorithm from Theorem~\ref{m1} and compute $\rk(G_i)$ for each $i$.

Then by Grushko's Theorem it follows that
$\rk(G)=\rk(G_1)+\dots+\rk(G_s)+r.$
\end{proof}

\section{Nielsen equivalence classes of generating tuples}

In this section we prove Theorem~\ref{B} from the Introduction.

\begin{thm}\label{thm:nielsen}
  Let $G$ be a group from the class $\mathcal M$. Then for each $k\ge
  1$ there exists only finitely many Nielsen-equivalence classes of
  $k$--tuples generating $G$.
\end{thm}
\begin{proof}
  First, observe that it suffices to prove Theorem~\ref{thm:nielsen}
  for groups from the class $\mathcal M_1$. Indeed, Grushko's theorem
  (see the proofs, for example, in \cite{KMW,W02}) ensures that every
  tuple generating a group $G_1\ast G_2$ is Nielsen equivalent to a
  tuple $(T_1,T_2)$ where $T_1$ generates $G_1$ and $T_2$ generates
  $G_2$. Note also that the statement of Theorem~\ref{thm:nielsen}
  holds for infinite cyclic groups (and, moreover, for finitely
  generated free groups).
  
  Thus assume that $G=\langle S|R\rangle\in \mathcal M_1$.
  Theorem~\ref{constC} implies that there is a constant $C>0$ such
  that if $T=(g_1,\dots, g_k)$ generates $G$ than either there exists
  a $k$--tuple $T'=(y_1,y_2,\dots, y_k)$ Nielsen-equivalent to $T$ such
  that $|y_i|_S\le C$ for $i=1,\dots, k$, or there exists a $k$--tuple
  $T'=(y_1,y_2,\dots, y_k)$ Nielsen-equivalent to $T$ such that for
  some $j<k$ the set $(y_1,\dots, y_j)$ generates a virtual fiber $U$
  of $G$ and $|y_i|_S\le C$ for $i=1,\dots, j$.
  
  In the former case we are done, so suppose that the latter occurs.
  Thus there is a subgroup $G_1$ of finite index in $G$ such that
  $G_1=H\rtimes \langle t\rangle$ and that the fiber $H$ of $G_1$ is
  commensurable with $U$. By Remark~\ref{rem:vf} we may assume that
  $G_1$ is normal in $G$, since a finite cover of a closed 3--manifold
  fibering over a circle also fibers over a circle and the two fiber
  groups are commensurable.
  
  Put $H_U:=H\cap U$. Note that $H_U$ is commensurable with both $H$
  and $U$.
  
  There is a subgroup $P\le H_U$ of finite index such that $P$ is
  normal in $G_1$. Then by Lemma~\ref{conjug} there are only finitely
  many conjugates of $H$ and of $P$ in $G$. Hence there are only
  finitely many conjugates of $H_U$ in $G$.
  
  Note also that since there are only finitely many subgroups $U$ of
  $G$ under consideration, only finitely many possibilities for $G_1$,
  $H$ and $H_U$ arise.
  
  We consider the following two cases:

{\bf Case 1}\qua The subgroup $H_U$ is commensurable with
$g^{-1}H_Ug$ for every $g\in G$. Then the group $\displaystyle
N=\cap_{g\in G} g^{-1}H_Ug$ is a subgroup of finite index in $H_U$ and
in $U$ that is normal in $G$. Therefore the group $G/N$ is virtually
cyclic.

The same argument as in the proof of the solvability of the relative
rank problem with respect to virtual fibers (see Proposition~\ref{vc}
and Remark~\ref{rem:vc}) now yields the conclusion that there are only
finitely many Nielsen-equivalence classes of $k$--tuples
$T'=(y_1,y_2,\dots, y_k)$ with these properties.

{\bf Case 2}\qua
Suppose that there is some $g\in G$ such that $H_U$ is not
commensurable with $gH_Ug^{-1}$. Recall that $U$ is generated by
$(y_1, \dots, y_j)$ and that $U$ is commensurable with $H_U$.  Thus
there is $q>j$ such that $H_U$ is not commensurable with
$y_q^{-1}H_Uy_q$.  Note that $y_q^{-1}H_Uy_q\le G_1$ since $G_1$ is
normal in $G$.  Recall also that $H_U$ has only finitely many
conjugates in $G$ and hence only finitely many possibilities for
$y_q^{-1}H_Uy_q$ arise.

The subgroup $H_U$ contains a subgroup $P$ of finite index in both $H$
and $H_U$ such that $P$ is normal in $G_1$. Hence by
Lemma~\ref{conjug} every subgroup containing $P$ is either
commensurable with $H_U$ or has finite index in $G$. By assumption
$y_q^{-1}H_Uy_q$ is not commensurable with $H_U$.  Hence $L:=\langle
H_U, y_q^{-1}H_Uy_q\rangle$ has finite index in $G_1$ and in $G$.
Indeed, otherwise $L$ is a finite extension of $H_U$ and $L$ contains
a virtual fiber $y_q^{-1}H_Uy_q$. Since virtual fibers are surface
groups, whenever one of them is contained in another, it must have
finite index in this bigger virtual fiber. Hence $y_q^{-1}H_Uy_q$ is a
subgroup of finite index in $L$. Since $H_U$ also has finite index in
$L$, this implies that $H_U$ is commensurable with $y_q^{-1}H_Uy_q$,
contrary to our assumptions. Thus indeed $L$ has finite index in
$G_1$.

Write $y_q$ as $y_q=g'ht^s$ where $g'$ is an element from a fixed
finite transversal for $G_1$ in $G$, where $h\in H$ and $s\in \mathbb
Z$. Since $L\le G_1$ has finite index in $G_1$, we conclude that
$y_q^{-1}H_Uy_q$ is not contained in $H$. Therefore there is $h_u\in
H_U$ such that $y_q^{-1}h_uy_q=h't^a$ where $h'\in H$ and $a\ne 0$.
Then
$$(g'h)^{-1} h_u (g'h)=t^{s} h' t^a t^{-s}=h''t^a \quad \text{ where }
h''\in H.$$
Hence for every integer $m$
$$(g'h)^{-1} h_u^m (g'h)=h_mt^{am} \quad \text{ where } h_m\in H,$$
that is
$$h_u^m g'h=g'hh_mt^{am}.$$
Put
$$y_q'=h_u^my_q=h_u^m g'ht^s=g'hh_mt^{am+s}=g't^{am+s}h_m'$$
where $h_m'\in H$.
We can choose an appropriate $m$ to ensure that $|am+s|\le |a|$.

Replacing $y_q$ by $y_q'=h_uy_q$ in $T'$ results in a Nielsen-equivalent
tuple since $h_u\in H_U\le U$. After this we can use the
subgroup of finite index $L$ in $G$ to make all the entries $y_i$ for
$i>j$, $i\ne q$ short. Next, using the fact that $H_U$ has finite index
in $H$, we multiply $h_q'=g't^{am+s}h_m'$ on the right by an
appropriate element of $H_U$ to make this entry short as
well.\end{proof}


\begin{thebibliography}{XXX}
  
\bibitem{A}
\textbf{I Agol},
  \emph{Tameness of hyperbolic 3--manifolds},
  \arxiv{math.GT/0405568}

\bibitem{ABC}
\textbf{J Alonso}, \textbf{T Brady}, \textbf{D Cooper}, \textbf{V Ferlini},
\textbf{M Lustig}, \textbf{M Mihalik}, \textbf{M Shapiro}, \textbf{H Short},
\emph{Notes on hyperbolic groups},
  from: ``Group theory from a geometric viewpoint (Trieste 1990)'',
  World Scientific, River Edge, NJ (1991) 3--63
  \MR{1170363}
  
\bibitem{Ar}
\textbf{G Arzhantseva},
\emph{On quasiconvex subgroups of word hyperbolic groups}
  Geom. Dedicata 87 (2001) 191--208
  \MR{1866849}
  
\bibitem{AO} \textbf{G Arzhantseva}, \textbf{A\,Yu Ol'shanskii},
\emph{Generality of the class of groups in which subgroups with a
lesser number of generators are free} Mat. Zametki 59 (1996) 489--496,
638 (Russian); translation in Math. Notes 59 (1996) 350--355
\MR{1445193}


\bibitem{BMS}
\textbf{G Baumslag}, \textbf{C\,F Miller III}, \textbf{H Short},
\emph{Unsolvable problems about small cancellation and word hyperbolic
  groups}
  Bull. London Math. Soc. 26 (1994) 97--101
  \MR{1246477}
  
\bibitem{BZ}
\textbf{M Boileau}, \textbf{H Zieschang},
\emph{Heegaard genus of closed orientable Seifert $3$--manifolds}
  Invent. Math. 76 (1984) 455--468
  \MR{0746538}

  
\bibitem{CG}
\textbf{D Calegari}, \textbf{D Gabai},
\emph{Shrinkwrapping and the taming of hyperbolic 3--manifolds},
\arxiv{math.GT/0407161}

  
\bibitem{C}
\textbf{R Canary},
\emph{A covering theorem for hyperbolic $3$--manifolds and its applications},
  Topology 35 (1996) 751--778
  \MR{1396777}

  
\bibitem{DT}
\textbf{N Dunfield}, \textbf{W\,P Thurston},
\emph{The virtual Haken conjecture: Experiments and examples},
  \gtref7{2003}{12}{399}{441}
  \MR{1988291}
  
\bibitem{E-T}
\textbf{D\,B\,A Epstein}, \textbf{J\,W Cannon}, \textbf{D\,F Holt},
\textbf{S\,V\,F Levy}, \textbf{M\,S Paterson}, \textbf{W\,P Thurston},
\emph{Word processing in groups},
  Jones and Bartlett Publishers, Boston, MA (1992)
  \MR{1161694}

  
\bibitem{Ga}
\textbf{D Gabai},
\emph{Convergence groups are Fuchsian groups},
  Ann. of Math. 136 (1992) 447--510
  \MR{1189862}

  
\bibitem{GH}
\textbf{\'E Ghys}, \textbf{P de la Harpe} (editors),
\emph{Sur les groupes hyperboliques d'apr\`es Mikhael Gromov},
  Progress in Mathematics 83,
  Birkh\"auser Boston Inc. Boston, MA, USA (1990)
  \MR{1086648}

  
\bibitem{Ger}
\textbf{V Gerasimov},
\emph{Detecting connectedness of the boundary of a hyperbolic group},
preprint (1999)
  
\bibitem{GS}
\textbf{S\,M Gersten}, \textbf{H Short},
\emph{Rational subgroups of biautomatic groups},
  Ann. of Math. 134 (1991) 125--158
  \MR{1114609}
  
\bibitem{G}
\textbf{M Gromov},
\emph{Hyperbolic groups},
  from: ``Essays in group theory'' (S\,M Gersten, editor),
  Math. Sci. Res. Inst. Publ. 8, Springer--Verlag (1987) 75--263
 \MR{0919829}
 
\bibitem{JR}
\textbf{W Jaco}, \textbf{H Rubinstein},
\emph{0--efficient triangulations of 3--manifolds},
  J. Differential Geom. 65 (2003) 61--168
  \MR{2057531}

\bibitem{J}
\textbf{K Johannson},
\emph{Topology and combinatorics of 3--manifolds},
  Lecture Notes in Mathematics 1599,
  Springer--Verlag, Berlin, Germany (1995)
  \MR{1439249}

\bibitem{Ka96}
\textbf{I Kapovich},
\emph{Detecting quasiconvexity: algorithmic aspects},
  from: ``Geometric and computational perspectives on infinite groups
  (Minneapolis, MN and New Brunswick, NJ, 1994)'',
  DIMACS Ser. Discrete Math.  Theoret. Comput. Sci. 25,
  Amer. Math. Soc., Providence, RI, USA (1996) 91--99
  \MR{1364182}

\bibitem{KMW} \textbf{I Kapovich}, \textbf{A Myasnikov}, \textbf{R
Weidmann}, \emph{Foldings, graphs of groups and the membership
problem}, Internat. J. Alg. Comput.  15 (2005) 95--128

\bibitem{KS}
\textbf{I Kapovich}, \textbf{H Short},
\emph{Greenberg's theorem for quasiconvex subgroups of word hyperbolic
  groups},
  Canad. J. Math. 48 (1996) 1224--1244
  \MR{1426902}

\bibitem{KW}
\textbf{I Kapovich}, \textbf{R Weidmann},
  \emph{Nielsen Methods for groups acting on hyperbolic spaces},
  Geom. Dedicata 98 (2003) 95--121
  \MR{1988426}

\bibitem{KW04}
\textbf{I Kapovich}, \textbf{R Weidmann},
\emph{Freely indecomposable groups acting on hyperbolic spaces},
  Internat. J. Algebra Comput. 14 (2004) 115--171
  \MR{2058318}

\bibitem{Ka01}
\textbf{M Kapovich},
\emph{Hyperbolic manifolds and discrete groups},
  Progress in Mathematics 183,
  Birkh\"auser, Boston, MA, USA (2001)
  \MR{1792613}


\bibitem{KMS}
\textbf{A Karrass}, \textbf{W Magnus}, \textbf{D Solitar},
\emph{Combinatorial Group Theory: Presentations of groups in terms of
  generators and relations},
  Interscience Publishers [John Wiley and Sons, Inc.],
  New York--London--Sydney (1966)
  \MR{0207802}


\bibitem{P}
\textbf{P Papasoglu},
\emph{An algorithm detecting hyperbolicity},
  from: ``Geometric and computational perspectives on infinite groups
  (Minneapolis, MN and New Brunswick, NJ 1994)'',
  DIMACS Ser. Discrete Math. Theor. Comput. Sci. 25,
  Amer. Math. Soc. Providence, RI, USA (1996) 193--200
  \MR{1364175}

\bibitem{R}
\textbf{E Rips},
\emph{Subgroups of small cancellation groups},
  Bull.  Lond. Math. Soc. 14 (1982) 45--47
  \MR{0642423}

\bibitem{Ru}
\textbf{J Rubinstein},
\emph{Polyhedral minimal surfaces, Heegaard splittings and decision
  problems for $3$--dimensional manifolds},
  from: ``Geometric topology (Athens, GA, 1993)'',
  AMS/IP Stud. Adv. Math. 2,
  Amer. Math. Soc. Providence, RI, USA (1997) 1--20
  \MR{1470718}

\bibitem{Sch}
\textbf{P Schupp},
\emph{Coxeter groups, 2--completion, perimeter reduction and subgroup
  separability},
  Geom. Dedicata 96 (2003) 179--198
  \MR{1956839}

\bibitem{Se}
\textbf{Z Sela},
\emph{The isomorphism problem for hyperbolic groups I},
  Ann. of Math. 141 (1995) 217--283
  \MR{1324134}

\bibitem{Si}
\textbf{C Sims},
\emph{Computation with finitely presented groups},
  volume 48 of \emph{Encyclopedia of Mathematics and its Applications},
  Cambridge University Press, Cambridge (1994)
  \MR{1267733}


\bibitem{So}
\textbf{J Souto},
\emph{The rank of the fundamental group of hyperbolic 3--manifolds
  fibering over the circle}, preprint (2005)


\bibitem{Tu}
\textbf{P Tukia},
\emph{Homeomorphic conjugates of Fuchsian groups},
  J. Reine Angew. Math. 391 (1988) 1--54
  \MR{0961162}


\bibitem{W02}
\textbf{R Weidmann},
\emph{The Nielsen method for groups acting on trees},
  Proc. London Math. Soc. 85 (2002) 93--118
\MR{1901370}

\bibitem{W04}
\textbf{R Weidmann},
\emph{The rank problem for sufficiently large Fuchsian groups},
preprint (2004)


\end{thebibliography}
\end{document}